\documentclass[12pt]{amsart}
\usepackage{amscd,amssymb,times,fullpage}
\newtheorem{theorem}{Theorem}

\newtheorem{proposition}[theorem]{Proposition}

\theoremstyle{definition}
\newtheorem{definition}[theorem]{Definition}

\newtheorem{example}[theorem]{Example}

\theoremstyle{remark}
\newtheorem{remark}[theorem]{Remark}

\newcommand\Ker{\operatorname{Ker}}

\catcode`\@=\active 

\begin{document}

\title{Contact pair structures and associated metrics}

\author{G.~Bande}
\address{Dipartimento di Matematica e Informatica, Universit{\`a} degli studi di Cagliari, Via Ospedale 72, 09124 Cagliari, Italy}
\email{gbande{\char'100}unica.it}
\author{A.~Hadjar}
\address{Laboratoire de Math{\'e}matiques, Informatique et
Applications, Universit{\'e} de Haute Alsace - 4 Rue de Fr{\`e}res
Lumi{\`e}re, 68093 Mulhouse Cedex, France}
\email{Amine.Hadjar{\char'100}uha.fr}
\thanks{The first author was supported by the MIUR Project: \textit{ Riemannian metrics and differentiable manifolds}--P.R.I.N. $2005$}


\begin{abstract}
We introduce the notion of contact pair structure and the
corresponding associated metrics,
 in the same spirit of the geometry of almost contact structures. We prove that, with respect to these
metrics, the integral curves of the Reeb vector fields are
geodesics and that the leaves of the Reeb action are totally
geodesic. Moreover, we show that, in the case of a metric contact
pair with decomposable endomorphism, the characteristic foliations
are orthogonal and their leaves carry induced contact metric
structures.
\end{abstract}

\maketitle

\section{Introduction}

A \textit{contact pair} on a smooth manifold $M$ is a pair of one-forms
$\alpha_1$ and $\alpha_2$ of constant and complementary classes,
for which $\alpha_1$ restricted to the leaves of
the characteristic foliation of $\alpha_2$ is a contact form and vice versa. This
definition was introduced in \cite{Bande1,BH} and is the analogous to that of
\textit{contact-symplectic pairs} and \textit{symplectic pairs} (see \cite{Bande2,BK,KM}).

In this paper we introduce the notion of \textit{contact pair structure}, that is a
contact pair on a manifold $M$ endowed with a tensor field $\phi$, of type $(1,1)$, such
 that $\phi^2=-Id + \alpha_1 \otimes Z_1 + \alpha_2 \otimes Z_2$ and $\phi(Z_1)=\phi(Z_2)=0$, where $Z_1$ and $Z_2$
are the Reeb vector fields of the pair.

We define \textit{compatible} and
\textit{associated} metrics and we prove several properties of these metrics.
For example, we show that the orbits of the locally free $\mathbb{R}^2$-action generated by
the two commuting \textit{Reeb vector fields} are totally geodesic with respect to these metrics.
Another important feature is that, given a metric which is associated with respect to a contact pair structure $(\alpha_1 , \alpha_2 , \phi)$,
if we assume that the endomorphism $\phi$ preserves the tangent spaces of the leaves of the two
characteristic foliations, then the contact forms induced on the
leaves are \textit{contact metric structures}.

In what follows we denote by $\Gamma(B)$ the space of sections of
a vector bundle $B$. For a given foliation $\mathcal{F}$ on a
manifold $M$, we denote by $T\mathcal{F}$ the subbundle of $TM$
whose fibers are given by the distribution tangent to the leaves.
All the differential objects considered are supposed smooth.

\section{Preliminaries on contact pairs}\label{s:prelim}

We gather in this section the notions concerning contact pairs that will
be needed in the sequel. We refer the reader
to \cite{Bande2, BH, BGK, BK} for further informations and
several examples of such structures.
\begin{definition}\label{d:cpair}
A pair $(\alpha_1, \alpha_2)$ of $1$-forms on a manifold is said
to be a contact pair of type $(h,k)$ if:
\begin{enumerate}
\item[i)] $\alpha_1\wedge
(d\alpha_1)^{h}\wedge\alpha_2\wedge (d\alpha_2)^{k}$ is a volume
form,
\item[ii)] $(d\alpha_1)^{h+1}=0$,
\item[iii)] $(d\alpha_2)^{k+1}=0$.
\end{enumerate}
\end{definition}
Since the form $\alpha_1$ (resp. $\alpha_2$) has constant class $2h+1$
(resp. $2k+1$), the distribution $\Ker \alpha_1 \cap \Ker d\alpha_1$ (resp.
 $\Ker \alpha_2 \cap \Ker d\alpha_2$) is completely integrable and then it determines the so-called characteristic
  foliation $\mathcal{F}_1$ (resp. $\mathcal{F}_2$) whose leaves are endowed with a contact form induced by $\alpha_2$ (resp. $\alpha_1$).\\

To a contact pair $(\alpha_1, \alpha_2)$ of type $(h,k)$ are
associated two commuting vector fields $Z_1$ and $Z_2$, called \textit{Reeb
vector fields} of the pair, which are uniquely determined by the
following equations:
\begin{eqnarray*}
&\alpha_1 (Z_1)=\alpha_2 (Z_2)=1  , \; \; \alpha_1
(Z_2)=\alpha_2
(Z_1)=0 \, , \\
&i_{Z_1} d\alpha_1 =i_{Z_1} d\alpha_2 =i_{Z_2}d\alpha_1=i_{Z_2}
d\alpha_2=0 \, ,
\end{eqnarray*}
where $i_X$ is the contraction with the vector field $X$.
In particular, since the Reeb vector fields commute, they determine a locally free
$\mathbb{R}^2$-action, called  the \textit{Reeb action}.

The kernel distribution of $d\alpha_1$ (resp. $d\alpha_2$) is
also integrable and then it defines a foliation whose leaves
inherit a contact pair of type $(0,k)$ (resp. $(h,0)$).

A contact pair of type $(h,k)$ has a local model (see \cite{BH, Bande2}), which means that for every point of the manifold there exists a coordinate
chart on which the pair can be written as follows:
$$
\alpha_1=dx_{2h+1}+\sum_{i=1}^{h}x_{2i-1}dx_{2i}  ,\; \;
\alpha_2=dy_{2k+1}+\sum_{i=1}^{k}y_{2i-1}dy_{2i} \, ,
$$
where $(x_1, \cdots ,x_{2h+1}, y_1 , \cdots , y_{2k+1})$ are the standard coordinates on $\mathbb{R}^{2h+2k+2}$.\\
The tangent bundle of a manifold $M$ endowed with a contact pair
$(\alpha_1, \alpha_2)$ can be split in different ways. For
$i=1,2$, let $T\mathcal F _i$ be the subbundle determined by the
characteristic foliation of $\alpha_i$, $T\mathcal G_i$ the
subbundle of $TM$ whose fibers are given by $\ker d\alpha_i \cap
\ker \alpha_1 \cap \ker \alpha_2$ and $\mathbb{R} Z_1, \mathbb{R}
Z_2$ the line bundles determined by the Reeb vector fields. Then
we have the following splittings:
$$
TM=T\mathcal F _1 \oplus T\mathcal F _2 =T\mathcal G_1 \oplus
T\mathcal G_2 \oplus \mathbb{R} Z_1 \oplus \mathbb{R} Z_2 .
$$
Moreover we have $T\mathcal F _1=T\mathcal G_1 \oplus \mathbb{R} Z_2 $ and $T\mathcal F _2=T\mathcal G_2 \oplus \mathbb{R} Z_1 $.

\section{Contact pair structures and almost contact
structures}\label{s2} In this section we firstly recall the basic
definitions of almost contact structures and their associated metrics. Next we introduce a
new structure, namely the \textit{contact pair structure}. More details
on almost contact structures can be found in
\cite{Blairbook}.
\subsection{Almost contact structures}
An almost contact structure on a manifold $M$ is a triple $(\alpha, Z, \phi)$ of a
one-form $\alpha$, a vector field $Z$ and a field of endomorphisms $\phi$ of the tangent bundle of $M$,
 such that $\phi ^2 =-Id + \alpha \otimes Z$, $\phi(Z)=0$ and $\alpha(Z)=1$. In particular, it follows that $\alpha \circ \phi=0$ and that
 the rank of $\phi$ is $\dim M -1$.

In the study of almost contact structures, there are two types of metrics which are particularly interesting and we recall their definition:
\begin{definition}
For a given almost contact structure $(\alpha, Z, \phi)$ on a manifold $M$, a Riemannian metric $g$ is called:
\begin{itemize}
\item[i)] compatible if $g(\phi X,\phi Y)=g(X,Y)-\alpha (X) \alpha(Y)$ for every $X,Y \in \Gamma (TM)$;
\item[ii)]associated if $g(X,\phi Y)=d\alpha (X,Y)$ and $g(X,Z)=\alpha (X)$ for every $X , Y \in \Gamma (TM)$.
\end{itemize}
\end{definition}
In particular if $g$ is an associated metric with respect to $(\alpha, Z, \phi)$, then $\alpha$ must be a contact form. For a given almost
contact structure there always exists a compatible metric and  for a given contact form there always exists an associated metric.
\begin{definition}
A contact metric structure $(\alpha, Z, \phi, g)$ is an almost contact structure $(\alpha, Z, \phi)$, where $\alpha$ is a contact form and $Z$ its Reeb
vector field, together with an associated metric $g$.
\end{definition}

\subsection{Contact pair structures} Now we want to generalize the notion of contact metric structure to the contact pairs.
To do that, we begin with the following definition:
\begin{definition}
A \emph{contact pair structure} on a manifold $M$ is a triple
$(\alpha_1 , \alpha_2 , \phi)$, were $(\alpha_1 , \alpha_2)$ is a
contact pair and $\phi$ a tensor field of type
$(1,1)$ such that:
\begin{gather}
\label{d:cpstructure}\phi^2=-Id + \alpha_1 \otimes Z_1 + \alpha_2 \otimes Z_2 ,\\
\label{d:cpstructure2}\phi(Z_1)=\phi(Z_2)=0 .
\end{gather}
where $Z_1$ and $Z_2$ are the Reeb vector fields of $(\alpha_1 ,
\alpha_2)$.
\end{definition}
As in the case of an almost contact structure, it is easy to check
that $\alpha_i \circ \phi =0$, $i=1,2$ and the rank of $\phi$ is
equal to $\dim M -2$.
\begin{remark}
Actually, the condition $\phi(Z)=0$ in the definition of almost
contact structure is not needed, since it follows from $\phi^2=-I
+\alpha \otimes Z$ (see \cite{Blairbook}, Theorem $4.1$ p. $33$).
In the case of a contact pair,
 it is possible to construct an endomorphism satisfying \eqref{d:cpstructure} but not \eqref{d:cpstructure2}. This can be done, for example,
 by taking $\phi$ to be an almost complex structure on $T\mathcal G_1 \oplus T\mathcal G_2$ and extending it by setting $\phi (Z_1)=- \phi (Z_2)= Z_1+Z_2$.
\end{remark}

On a manifold $M$ endowed with a contact pair, there always exists
an endomorphism $\phi$ verifying \eqref{d:cpstructure} and
\eqref{d:cpstructure2}, since on the subbundle $T\mathcal G_1
\oplus T\mathcal G_2$ of $TM$, the $2$-form $d\alpha_1 +
d\alpha_2$ is symplectic. Then one can choose an almost complex
structure on $T\mathcal G_1 \oplus T\mathcal G_2$ and extend it to
be zero on the Reeb vector fields, to produce an endomorphism
$\phi$ of $TM$ with the desired properties.

The following definition is justified by the fact that we are interested on the contact structures induced on the leaves of
the characteristic foliations of a contact pair.
\begin{definition}
The endomorphism $\phi$ is said to be \textit{decomposable} if
$\phi (T\mathcal{F}_i) \subset T\mathcal{F}_i$, for $i=1,2$.
\end{definition}
The condition for $\phi$ to be decomposable is equivalent to $\phi (T\mathcal{G}_i)= T\mathcal{G}_i$, because
$\phi$ vanishes on $Z_1$, $Z_2$ and has rank $2k$ (resp. $2h$) on $T\mathcal{F}_1$ (resp. $T\mathcal{F}_2$).

By Definition \ref{d:cpair}, $\alpha_1$ induces a contact form on the leaves of the
characteristic foliation of $\alpha_2$ and vice versa. Then it is easy to prove the following:
\begin{proposition}
If $\phi$ is decomposable, then $(\alpha_1 , Z_1 ,\phi)$ (resp. $(\alpha_2 , Z_2 ,\phi)$) induces an almost contact structure
on the leaves of $\mathcal{F}_2$ (resp.$\mathcal{F}_1$).
\end{proposition}
Here is a simple example of contact pair on a manifold which is just a product of two contact manifolds and
$\phi$ is not decomposable:
\begin{example}\label{ex1}
Let us consider standard coordinates $(x_1, y_1, x_2, y_2,z_1, z_2
)$ on $\mathbb{R} ^6$ and the contact pair $(\alpha_1 , \alpha_2)$
given by:
$$
\alpha_1= dz_1 -x_1 dy_1 \, \, \, \, \, \, \alpha_2= dz_2 -x_2
dy_2 ,
$$
with Reeb vector fields $Z_1=\frac{\partial}{\partial z_1}$
, $Z_2=\frac{\partial}{\partial z_2}$. Let us define $\phi$ as follows:
\begin{eqnarray*}
&\phi \left(a_1 \frac{\partial}{\partial x_1}+ b_1
\frac{\partial}{\partial y_1}+  a_2 \frac{\partial}{\partial x_2}+
b_2\frac{\partial}{\partial y_2}+ c_1 \frac{\partial}{\partial
z_1}+ c_2\frac{\partial}{\partial z_2} \right)=\\
& = \left(-a_2 \frac{\partial}{\partial x_1}-b_2
\frac{\partial}{\partial y_1}+ a_1 \frac{\partial}{\partial x_2}+
b_1\frac{\partial}{\partial y_2}-x_1 b_2 \frac{\partial}{\partial
z_1}+ x_2 b_1\frac{\partial}{\partial z_2} \right) \, ,
\end{eqnarray*}
where $a_1 , a_2, b_1,b_2 , c_1, c_2 $ are smooth functions on $\mathbb{R}^6$.\\
Then $(\alpha_1, \alpha_2, \phi)$ is a contact pair structure and
$\phi$ is not decomposable. This because, for example,
$\frac{\partial}{\partial x_1} \in \Ker (\alpha_2 \wedge d\alpha_2) \, , \, \frac{\partial}{\partial x_2} \in \Ker (\alpha_1 \wedge d\alpha_1)$ but
 $\phi (\frac{\partial}{\partial x_1})=\frac{\partial}{\partial x_2}$.
\end{example}
\begin{remark}\label{generalstructure}
A more general structure, similar to the almost contact structures,
is obtained by considering a 5-tuple $(\alpha_1, \alpha_2 ,
Z_1 , Z_2 , \phi)$, where $\phi^2 = -Id + \alpha_1 \otimes Z_1 + \alpha_2
\otimes Z_2$, $ \phi(Z_1)=\phi(Z_2)=0$ and the $\alpha_i$'s are just non-vanishing $1$-forms
verifying $\alpha_i (Z_j)= \delta_{ij}$.
\end{remark}
\section{Compatible and associated metrics}\label{s3}

For a given contact pair structure on a manifold $M$, it is
natural to consider the following metrics:
\begin{definition}
Let $(\alpha_1 , \alpha_2 ,\phi )$ be a contact pair structure on a manifold $M$, with Reeb vector fields $Z_1$ and $Z_2$.
A Riemannian metric $g$ on $M$ is called:
\begin{itemize}
\item[i)] compatible if $g(\phi X,\phi Y)=g(X,Y)-\alpha_1 (X) \alpha_1 (Y)-\alpha_2 (X)
\alpha_2 (Y)$ for all $X,Y \in \Gamma (TM)$,

\item[ii)] associated if $g(X, \phi Y)= (d \alpha_1 + d \alpha_2) (X,Y)$ and $g(X,
Z_i)=\alpha_i(X)$, for $i=1,2$ and for all $X,Y \in \Gamma (TM)$. \label{ass-metric}
\end{itemize}

\end{definition}
It is clear that an associated metric is also compatible, but the
converse is not true. In the Example \ref{ex1}, the metric $g=
\sum _{i=1} ^2 (dx_i ^2 +dy_i ^2+\alpha_i ^2)$ is compatible but
not associated.
\begin{definition}
A \emph{metric contact pair} (MCP) on a manifold $M$ is a $4$-tuple $(\alpha_1, \alpha_2, \phi, g)$ where $(\alpha_1, \alpha_2, \phi)$ is
a contact pair structure and $g$ an associated metric with respect to it.
\end{definition}

Like for compatible metrics on almost contact manifolds, we have:
\begin{proposition}
For every contact pair structure on a manifold $M$ there exists a compatible metric.
\end{proposition}
\begin{proof}
Let $(\alpha_1 , \alpha_2, \phi)$ be a contact pair structure on $M$.
Pick any Riemannian metric $h$ on $M$ and define $k$ as
$$
k(X,Y)= h(\phi^2 X , \phi^2 Y) + \alpha_1 (X) \alpha_1
(Y)+\alpha_2 (X) \alpha_2 (Y) \, .
$$
Now put
$$
g(X,Y)=\frac{1}{2} (k(X,Y)+ k(\phi X, \phi Y) +\alpha_1 (X)
\alpha_1 (Y)+\alpha_2 (X) \alpha_2 (Y)) \, .
$$
It is straightforward to show that $g$ is a compatible metric with respect to
$(\alpha_1 , \alpha_2, \phi)$.
\end{proof}
Compatible and associated metrics have several interesting
properties given by the following
\begin{theorem}
Let $M$ be a manifold endowed with a contact pair structure $(\alpha_1 , \alpha_2, \phi)$, with Reeb vector fields $Z_1 , Z_2$. Let $g$ be
 a compatible metric with respect to it,
 with Levi-Civit\`a connection $\nabla$. Then we have:
\begin{enumerate}
\item $g(Z_i, X)= \alpha_i (X)$, for $i=1,2$ and for every $X \in
\Gamma(TM)$;

\item $g(Z_i, Z_j)= \delta_i ^j$, $i,j=1,2$;

\item $\nabla _{Z_i} Z_j = 0$, $i,j=1,2$ (in particular the integral curves of the
Reeb vector fields are geodesics);

\item the Reeb action is totally geodesic (i.e the orbits are totally geodesic $2$-dimensional submanifolds).

\end{enumerate}
Moreover, if $g$ is an associated metric and $L_{Z_i}$ is the Lie derivative along $Z_i$, then $L_{Z_i}\phi=0$ if
and only if $Z_i$ is a Killing vector field.
\end{theorem}
\begin{proof}
The first two properties are easy consequences of the definitions
of compatible metric and of the Reeb vector fields of a contact
pair. For the third property, let us remember that $\alpha_i$ is
invariant by $Z_j$, $i,j=1,2$, and then:
\begin{alignat*}{1}
0 &=(L_{Z_j} \alpha_i)(X)=Z_j(\alpha_i (X))- \alpha_i
(L_{Z_j} X)\\
&=Z_j(g(X, Z_i))- g (Z_i ,\nabla_{Z_j} X - \nabla _ X
Z_j)\\
&=g(\nabla_{Z_j}X, Z_i)+ g(X, \nabla_{Z_j}Z_i)- g (Z_i
,\nabla_{Z_j} X - \nabla _ X Z_j)\\
&=g(X, \nabla_{Z_j}Z_i)+g (Z_i
,\nabla _ X Z_j) \, .
\end{alignat*}
Summing up with $(L_{Z_i} \alpha_j)(X)$ and recalling that $[Z_i
, Z_j]=0$, we get $g(X, \nabla_{Z_j}Z_i)=0$ for all $X$.

To prove the fourth property, let us consider the
second fundamental form $B$ of an orbit $\tilde M$ of the Reeb action. Since the Reeb vector fields are tangent to the orbits of the Reeb action,
we can choose $\{Z_1 , Z_2\}$ (restricted to $\tilde M$) as a basis of the tangent space at every point of $\tilde M$ . For a vector field $X$
on $M$ along $\tilde M$,
 let us denote by $X^T$ its tangential component. Then, by the third property, we have:
$$
B(Z_i , Z_j)= \nabla_{Z_i} Z_j - (\nabla_{Z_i} Z_j )^T=0 \, .
$$
Finally, if $g$ is associated, we want to prove that $[L_{Z_i}g](X,Y )=0$ for all $X, Y \in \Gamma(TM)$ if and only if $(L_{Z_i}\phi)(Y)=0$
for every $Y \in \Gamma(TM)$.
 This is clearly true for $Y=Z_1$ or $Z_2$, since in this case $[L_{Z_i}\phi](Z_j )$ vanishes identically and
 $L_{Z_i}\alpha_j=0$ applied to a vector field $X$ implies $[L_{Z_i}g](X,Z_j )=0$.\\
It remains to prove the assertion for $Y$ in $\Ker \alpha_1 \cap \Ker \alpha_2$. In this case, for all vector fields $X$, $Y$, we
have:
\begin{alignat*}{1}
0&= L_{Z_i} (d\alpha_1+d\alpha_2)(X,Y)\\
&=Z_i
g(X,\phi Y)-(d\alpha_1+d\alpha_2)(L_{Z_i}X,Y)-(d\alpha_1+d\alpha_2)(X,L_{Z_i}Y)\\
&=[L_{Z_i}g](X,\phi Y )+g(X, [L_{Z_i}\phi](Y))
\end{alignat*}
and this completes the proof since $\phi$ restricted to $\Ker \alpha_1 \cap \Ker \alpha_2$ is an isomorphism.
\end{proof}

\begin{remark}
Codimension two geodesible foliations of closed $4$-manifolds have
been studied by Cairns and Ghys in \cite{ghys}. In particular,
they have proven that, in this situation, there exists a metric on the manifold for which
the leaves of the foliation are minimal and have same
constant curvature $K=0, 1, -1$. Our case belongs to what they have
called \textit{parabolic case}, that is $K=0$.
\end{remark}
As for the metric contact structures (see \cite{Blairbook} for
example), by polarization one can show the existence of associated
metrics, and in fact we have:
\begin{proposition}
For every contact pair $(\alpha_1, \alpha_2)$ on a manifold $M$, there exists an endomorphism $\phi$ of $TM$
and a metric $g$ such that $(\alpha_1, \alpha_2 , \phi, g )$ is a
metric contact pair. Moreover $\phi$ can be chosen to be
decomposable.
\end{proposition}
\begin{proof}
Take any Riemannian metric $k$. Since $d\alpha_1 + d\alpha_2$ is
symplectic on the subbundle $T\mathcal G_1 \oplus T\mathcal G_2$,
then it can be polarized by using $k$ to obtain a metric $\tilde g$ and
an almost complex structure $\tilde \phi$ compatible with $\tilde g$. Extending $\tilde \phi$ to be zero on the Reeb vector fields, we obtain
 $\phi$. Defining $g$ to be equal to $\tilde g$ on $T\mathcal G_1 \oplus T\mathcal G_2$ and putting
 $g(X, Z_i)=\alpha_i (X)$, gives the desired metric.

To obtain a decomposable $\phi$, polarize $d\alpha_i$ on
$T\mathcal G_i$  to obtain two metrics $\tilde g_i$ and two endomorphisms $\tilde \phi_i$ on $T\mathcal G_i$ ($i=1,2$)
and then take the direct sums. Finally, extend them as before to obtain
the desired tensors.
\end{proof}

We end this section with two results concerning the structures induced on the leaves of the characteristic foliations:
\begin{theorem}
Let $M$ be a manifold endowed with a MCP $(\alpha_1, \alpha_2 ,
\phi, g )$ and suppose that $\phi$ is decomposable. Then
$(\alpha_i, \phi , g)$ induces a contact metric structure on the
leaves of the characteristic foliation of $\alpha_j$ for $i \neq
j$, $i,j=1,2$.
\end{theorem}
\begin{proof}
We will prove the assertion only for $(\alpha_1, \phi , g)$, the other case being completely similar.
If $F$ is a leaf of the foliation $\mathcal{F}_2$,
 then $(\alpha_1 , \phi)$ induces an almost contact structure on it and $\alpha_1$ is a contact form. Since $g$ is an associated metric
 with respect to $(\alpha_1, \alpha_2 , \phi )$,
 by Definition \ref{ass-metric}, when restricted to vectors $X,Y$ which are tangent to $F$ we have:
$$
g(X, \phi Y)= d \alpha_1  (X,Y) \, \, , \, \, g(X,
Z_1)=\alpha_1(X) \, ,
$$
showing that its restriction to $F$ is an associated metric with respect to the contact form induced by $\alpha_1$.
\end{proof}
In a similar way, recalling that the characteristic foliation of $d\alpha_i$ for $i=1,2$ is given by $\ker d\alpha_i$, we can prove:
\begin{theorem}
Let $M$ be a manifold endowed with a MCP $(\alpha_1, \alpha_2 ,
\phi, g )$ and suppose that $\phi$ is decomposable. Then
$(\alpha_1, \alpha_2 , \phi , g)$ induces a metric contact pairs
of type $(h,0)$ on the leaves of the characteristic foliation of
$d\alpha_2$
 and one of type $(0,k)$ on the leaves of the characteristic foliation of $d\alpha_1$.
\end{theorem}

\subsection{Orthogonal foliations}
The following theorem explains the link between decomposability of
$\phi$ and orthogonality of the characteristic foliations when the
metric is an associated one:
\begin{theorem}
For a MCP $(\alpha_1 , \alpha_2, \phi, g)$, the tensor $\phi$ is
decomposable if and only if the foliations $\mathcal{F}_1 ,
\mathcal{F}_2$ are orthogonal.
\end{theorem}
\begin{proof}
Suppose that $\phi (T\mathcal{F}_i) \subset T\mathcal{F}_i$ for $i=1,2$ and let $X
\in \Gamma (T\mathcal{F}_1 ), Y \in
\Gamma(T\mathcal{F}_2)$. Because $g$ is associated, by the choice of
$X$, $Y$, we have:
$$
g(X, \phi Y )= d\alpha_1 (X,Y) + d\alpha_2 (X,Y)=0 \, \, , \, \, g(X,
Z_1)= \alpha_1 (X)=0 ,
$$
which proves that $X$ is orthogonal to $T\mathcal{F}_2$ and then the two foliations are orthogonal.

Conversely, suppose that the two foliations are orthogonal. Then, for
$X \in \Gamma (T\mathcal{F}_1)$ and every $Y \in \Gamma (T\mathcal{F}_2)$, we have
$$
g( Y, \phi X)=(d\alpha_1 +d\alpha_2) (Y,X)=0
$$
which implies that $\phi X \in \Gamma (T\mathcal{F}_2^{\bot})= \Gamma (T\mathcal{F}_1)$,
that is $\phi$ is decomposable.

\end{proof}

In the Example \ref{ex1} the foliations are orthogonal with respect to the metric $g=
\sum _{i=1} ^2 (dx_i ^2 +dy_i ^2+ \alpha_i ^2)$ which is
compatible but not associated because $\phi$ is not decomposable.

Here is an example of MCP with decomposable $\phi$ on a nilpotent Lie group and its corresponding nilmanifolds:

\begin{example}
Consider the simply connected nilpotent Lie group $G$ with structure
equations:
\begin{eqnarray*}
&d\omega_1= d\omega_6=0 \; \; , \; \; d\omega_2= \omega_5 \wedge
\omega _6 ,\\
&d\omega_3=\omega_1 \wedge \omega_4 \; \; , \; \;
d\omega_4= \omega_1 \wedge \omega_5 \; \; , \; \; d\omega_5 =
\omega_1 \wedge \omega_6 \, ,
\end{eqnarray*}
where the $\omega_i$'s form a basis for the cotangent space of $G$
at the identity.

The pair $(\omega_2 , \omega_3)$ is a contact pair of type $(1,1)$
with Reeb vector fields $(X_2, X_3)$, the $X_i$'s being dual to
the $\omega_i$'s. Now define $\phi$ to be zero on the Reeb vector
fields and
$$
\phi(X_5)=-X_6 \; \; , \; \; \phi(X_6)=X_5 \; \; , \; \; \phi
(X_1)= -X_4 \; \; , \; \; \phi (X_4)=X_1 \, .
$$
Then $\phi$ is easy verified to be decomposable and the metric $g=\sum_{i=1}^6 \omega_i ^2$ is an associated metric with respect to
$(\omega_2 , \omega_3, \phi)$.

Since the structure constants of the group are rational, there
exist lattices $\Gamma$ such that $G/ \Gamma$ is compact. Since
the MCP on $G$ is left invariant, it descends to all quotients $G/ \Gamma$ and we obtain nilmanifolds carrying the same type of structure.
\end{example}

%

\section*{Final comments}
Further metric properties of the contact pairs structures are studied in \cite{BHinprogress},
where we prove, for example, that the characteristic foliations are minimal with respect to an associated metric.

In \cite{BH2} we study the analog for a contact pair structure of the notion of normality for almost contact structures.
We give there several constructions involving Boothby-Wang fibrations and flat bundles, which can be used to produce
more examples of metric contact pairs.

\end{document}